\documentclass[]{amsart}
\usepackage{graphicx}
\usepackage{amsfonts}

\newtheorem{theorem}{Theorem}[section]
\newtheorem{lemma}[theorem]{Lemma}

\newtheorem{corollary}[theorem]{Corollary}
\newtheorem{proposition}{Proposition}[section]

\theoremstyle{remark}
\newtheorem{remark}[theorem]{Remark}

\theoremstyle{definition}
\newtheorem{definition}[theorem]{Definition}
\newtheorem{example}[theorem]{Example}

\newcommand{\an}{\mathcal{A}} 
\newcommand{\Ad}{\textup{Ad}}
\newcommand{\mcc}{M_{\bf C}}

\newcommand{\khat}{\widehat{K}}
\newcommand{\uhat}{\widehat{U}}

 \newcommand{\real}{{\bf R}}

\newcommand{\hh}{\mathcal{H}}

\newcommand{\cc}{{\bf C}}

\newcommand{\bbar}{\left[ \begin{array}}
\newcommand{\ebar}{\end{array} \right] }
\newcommand{\bdm}{\begin{displaymath}}
\newcommand{\edm}{\end{displaymath}}
\newcommand{\beq}{\begin{equation}}
\newcommand{\beqa}{\begin{eqnarray}}
\newcommand{\beqas}{\begin{eqnarray*}}
\newcommand{\eeq}{\end{equation}}
\newcommand{\eeqa}{\end{eqnarray}}
\newcommand{\eeqas}{\end{eqnarray*}}
\newcommand{\dd}{\textup{d}}

\begin{document}

\title[Curved flats, pluriharmonic maps and pseudo space forms]
{Curved flats, pluriharmonic maps and constant curvature immersions into pseudo-Riemannian space forms}


\author{David Brander}
\address{Department of Mathematics\\ Faculty of Science\\Kobe University\\1-1, Rokkodai, Nada-ku, Kobe 657-8501\\ Japan}
\email{brander@math.kobe-u.ac.jp}

\begin{abstract}
We study two aspects of the loop group formulation for isometric
immersions with flat normal bundle  of space forms.
The first aspect is to examine the loop group maps
 along different ranges of the loop
parameter.  This leads to various equivalences between global
isometric immersion problems among different space forms and
pseudo-Riemannian space forms.  As a corollary, we obtain a non-immersibility
theorem for spheres into certain pseudo-Riemannian spheres and
hyperbolic spaces.

  The second aspect pursued is to clarify the relationship between
 the loop group formulation of isometric immersions of space forms
and that of pluriharmonic maps into symmetric spaces.
We show that the objects in the first class
 are, in the real analytic case,
extended pluriharmonic maps into certain symmetric spaces
which satisfy an extra reality condition along a totally
real submanifold.  We show how to construct such pluriharmonic
maps for general symmetric spaces from curved flats, using
a generalised DPW method.

\end{abstract}
\keywords{Isometric immersions, space forms, pluriharmonic maps, loop groups}
\subjclass[2000]{Primary 37K10, 37K25, 53C42, 53B25; Secondary 53C35}

\maketitle
\section{Introduction} \label{intro}
It is well known that harmonic maps from a Riemann surface 
into a symmetric space are 
integrable systems which can be approached successfully using
loop group techniques, which followed from the work of Uhlenbeck
\cite{uhlenbeck1989} on harmonic maps from $S^2$ into a Lie group.
  Further, various geometrical problems,
such as constant mean curvature surfaces, have been studied
successfully by showing that they are such harmonic maps (see,
for example, \cite{fordywood} for an introduction).
In higher dimensions, pluriharmonic maps into symmetric spaces
were also shown to have a similar approach by
 Ohnita and Valli \cite {ohnitavalli}, although the 
applications to special submanifolds appear to be little explored
thus far.  It turns out, as we will show, that the loop group maps
corresponding to isometric immersions of space forms, which were
defined by Ferus and Pedit \cite{feruspedit1996}, are a special
case of such pluriharmonic maps.

\subsection{Isometric immersions of space forms}
In this paper, we first study, in Section \ref{part2},
the interpretations of these loop group maps
 for different ranges of the spectral parameter. 
We  show, in Theorem \ref{corprop},
that the  map corresponding to each isometric immersion
problem actually contains three families of immersions, into
different space forms and pseudo-Riemannian space forms, for
values of the parameter  along $\real^*$, $i\real^*$ or $S^1$,
denoting the non-zero real, imaginary and unitary numbers respectively.  
It is also observed (Remark \ref{other}) that 
constant curvature immersions with flat normal bundle into
other pseudo-Riemannian space forms, beyond those arising here,
have an analogous loop group formulation, and solutions
can be generated via the AKS theory, as in \cite{feruspedit1996},
or by the well known dressing procedure.

Constant curvature immersions with flat normal bundle 
into pseudo-Riemannian space forms  have previously been 
studied by Barbosa, Ferreira and Tenenblat \cite{tenetal}, and
Dafeng, Qing and Yi \cite{qdy}, with
some additional assumptions to ensure the existence of special
coordinates. 
 The loop group formulation given here has
what is perhaps an 
advantage, in that it is coordinate-free, and therefore applies to
all cases.

This formulation  is global, and we 
also prove, in Proposition \ref{metriclemma},
that, for a given loop group map, completeness is equivalent
among all the isometric immersions obtained from it.

A corollary of our results is the equivalence among
various global isometric immersion problems of space forms
into pseudo-Riemannian space forms, stated in Corollary \ref{apcor}.
Among the applications of this is a coordinate-free proof,
 Corollary \ref{apcor2}, of the known result that
the problems of globally isometrically immersing the hyperbolic space
$H^m$ into a space form $Q_{\tilde{c}}^{2m-1}$, for $-1 < \tilde{c}$,
with $\tilde{c} \neq 0$, are equivalent, that is, independent of $\tilde{c}$.
This is related to the conjecture that no such isometric immersion is
possible for any $\tilde{c} > -1$, an extension of Hilbert's result
regarding immersions of the hyperbolic plane into Euclidean 3-space.
Following work by Pedit and Xavier,
\cite{pedit}, \cite{xavier}, this conjecture has been proven under the assumption
that the fundamental group of  the manifold being immersed is non-trivial, by Nikolayevsky \cite{nikolayevsky}, but remains a conjecture for the simply connected case.
One can alternatively show that the cases are equivalent 
for all $\tilde{c} > -1$ by using JD Moore's global principle coordinates
\cite{moore1}, and reducing the problem to finding
a global solution for 
the generalised sine-Gordon equation, studied by Terng and 
Tenenblat in \cite{terngtenenblat}, \cite{terng1980}.

 Another application,
Corollary \ref{apcor3}, is that there is no global isometric immersion
with flat normal bundle of a sphere $S^m(R)$, of radius $R<1$, into
the pseudo-Riemannian sphere $S^{2m-1}_{m-1}$, or of a sphere of any
radius into the pseudo-Riemannian hyperbolic space $H_{m-1}^{2m-1}$.

\subsection{Relation to pluriharmonic maps}
The remainder of the paper, beginning in Section \ref{dpw},
is devoted to exploring the relations between the loop
group formulation for isometric immersions of space forms,
and that of pluriharmonic maps into a symmetric space.

We show, in Theorem \ref{theorem1}, 
how to construct special pluriharmonic maps whose
extended families satisfy an extra reality condition when
restricted to a certain totally real submanifold.
They are constructed from real analytic 
curved flats into a different symmetric space. Many examples of such
curved flats can be constructed via the AKS theory \cite{feruspedit1996II}. 

As a special case, we show, in Theorem \ref{maintheorem},
that the loop group maps for isometric immersions studied
in Section \ref{part2}
 are, in the real analytic case,  just
restrictions to a totally real submanifold of pluriharmonic
maps into certain symmetric spaces. 

The proof  given here uses the generalised DPW method
\cite{dorfmeisterpeditwu}, \cite{branderdorf}, which
associates a certain type of  loop group map, which includes
both the isometric immersions of space forms and pluriharmonic
maps, to a simpler map, essentially a curved flat.  This is
then extended to a holomorphic map from a complex manifold into
the loop group associated to the complexification of the Lie group.
 We then apply the DPW correspondence again,
and obtain a pluriharmonic map.
In general the DPW method can only be applied to elements
which are in the big cell of the loop group, 
but we are able to get around this problem in
this case by renormalizing at different points of $M$.

 This example suggests that pluriharmonic
maps into symmetric spaces, perhaps satisfying extra conditions such as an
additional reality condition, should yield solutions to other interesting 
problems in geometry.  

In Section \ref{application} we recharacterize the problem of identifying
which totally geodesic submanifolds with flat normal bundle of 
the sphere and hyperbolic space (which is a degenerate case in
the formulation of Ferus and Pedit)
 belong to the families defined in
\cite{feruspedit1996}, in terms of pluriharmonic maps.

%*************************************************

\section{Limited connection order maps into loop groups} \label{lc}
In this section we outline some definitions and terminology. For
further details, we refer the reader to \cite{branderdorf}.
Let $G$ be a complex Lie group, with 
Lie algebra $\mathfrak{g}$. Let 
$\Lambda G$ be the group of real analytic  maps from the 
unit circle $S^1$ into $G$, with a topology that makes $\Lambda G$
a Banach Lie group. Any element $\gamma$ of $\Lambda G$ has an extension to 
a holomorphic map into $G$ on some annulus, $\an_\gamma$, containing $S^1$, and
in fact all the examples considered here are holomorphic on $\cc \setminus \{0\} := \cc^*$. 
  Let $M$ be a smooth manifold 
 and denote by $\Lambda G(M)$ the group of smooth maps
 $M \to \Lambda G$ normalised
to the identity at some fixed base point $p \in M$.   
If $F \in \Lambda G(M)$, then  for each value of the loop parameter $\lambda \in \an_F$,
we can expand the Maurer-Cartan form 
$F_\lambda^{-1} \dd F_\lambda$  as a
Fourier series in $\lambda$,
\bdm
F_\lambda^{-1} \dd F_\lambda = \sum _i A_i \lambda ^i.
\edm
For any subgroup $\mathcal{H}$  of the loop group $\Lambda G$,  
and any extended integers $a\leq b$,
we define $\mathcal{H} (M)_a^b$
 to be the set of elements in $\mathcal{H} (M)$ whose
Maurer-Cartan form is of bottom and top degree $a$ and $b$ respectively in
$\lambda$,
and call elements of
these sets \emph{connection order $(_a^b)$ maps}.

Let $\mathcal{H}^0$ denote the subgroup of constant loops in
$\mathcal{H}$. Note that $\mathcal{H}^0$ is a subgroup of 
$G = \Lambda G^0$.  The natural objects of study in this paper
are maps into the quotient space $\frac{\mathcal{H}}{\mathcal{H}^0}$
(where the action is right multiplication),
and we denote the set of connection order $(_a^b)$ maps from $M$
into this group by  $\frac{\mathcal{H}}{\mathcal{H}^0}(M)_a^b$.
For a fixed value of $\lambda$ these are maps from
$M$ into some 
quotient space  $\frac{U}{\mathcal{H}^0}$, where $U$ is some
subgroup of $G$.

\subsection{Some further terminology}
Let $\Lambda ^\pm G$ denote the subgroups of $\Lambda G$ consisting
of loops which extend analytically  to $D^\pm$, where $D^+$ is the
unit disc and $D^-$ is the complement of its closure in the Riemann sphere.
$\Lambda ^+_e G$ is the subgroup of $\Lambda ^+G$ whose elements 
are normalised to the identity at $\lambda = 0$, and $\Lambda^-_e G$
is the analogue at $\lambda = \infty$.  If $\mathcal{H}$ is a 
subgroup of $\Lambda G$ then the intersection of $\mathcal{H}$
with these subgroups are accordingly denoted $\mathcal{H}^\pm$ and 
$\mathcal{H}^\pm_e$.
The loop group $\Lambda G$ is called \emph{Birkhoff decomposable}
if, on an open dense neighbourhood of the identity, called the big
cell $B \Lambda G$, there is a (left Birkhoff)
decomposition
\beq \label{decomp}
\Lambda G = (\Lambda ^+_e G) \cdot (\Lambda^-G),
\eeq
and the map from $\Lambda ^+_e G  \times \Lambda^-G$ to this open
set is an analytic diffeomorphism. There is also an analogous
right Birkhoff decomposition, substituting $\pm$ for $\mp$.
In this paper, the group $G$ is always complex semisimple, so $\Lambda G$
is Birkhoff decomposable \cite{pressleysegal}.  A subgroup $\mathcal{H}$ of $\Lambda G$
is called Birkhoff decomposable if, in the decomposition given by
(\ref{decomp}), both factors on the right hand side are also in
the subgroup.  

Let $\phi$ be an automorphism of a subgroup $\mathcal{H}$, which is an
extension  of an automorphism of
$\mathcal{H}^0$.  We will say that
$\phi$ is \emph{positive} or \emph{negative (holomorphic)} if it
takes $\mathcal{H}^\pm$ to $\mathcal{H}^\pm$
or $\mathcal{H}^\pm$ to $\mathcal{H}^\mp$ respectively.
Note that a negative automorphism must be of even order.
A commonly used positive automorphism is the twisting:
\beq  \label{sigma}
(\sigma X)(\lambda) := \sigma^0 X(-\lambda),
\eeq
where $\sigma^0$ is an involution of the Lie group $G$.

A positive or negative holomorphic
\emph{reality condition} on $\mathcal{H}$ is a positive or negative
holomorphic  extension to $\mathcal{H}$ of a
 reality condition on $\mathcal{H}^0$.
Examples we will use are: for $\rho^0$ and $\tau^0$  
 reality conditions on $G$, define positive $\rho$
and negative $\tau$ on $\Lambda G$ by 
\beqa \label{rho}
\rho x(\lambda) := \rho^0 \, x(\bar{\lambda}),\\
\tau x(\lambda) := \tau^0 \, x(1/\bar{\lambda}). \label{tau}
\eeqa
Note that $x$ is understood to be defined on some annulus $\an_x$ around $S^1$, and
so $x(1/\bar{\lambda}) \neq x(\lambda)$.  
For any automorphism $\phi$ of $\Lambda G$ denote its fixed point 
subgroup by $\Lambda G _\phi$. Elements of  
$\Lambda G_\rho$ are real (that is take values 
in the real form  of $G$ determined by $\rho^0$) for real
values of $\lambda \in \an$, while elements of 
$\Lambda G_\tau$ are real for $\lambda \in S^1$.
For every  loop group automorphism
 used in
this paper, the corresponding automorphism of the Lie algebra
is given by the same formula, so we will use the same notation
to denote it. 

The important fact to note is that for any Birkhoff
decomposable subgroup $\mathcal{H}$ and any \emph{positive} 
holomorphic finite order automorphism $\phi$, the fixed point subgroup
$\mathcal{H}_ \phi$ is also Birkhoff decomposable. The same
does not hold if $\phi$ is
negative, as, given a Birkhoff decomposition $x = x_+y_-$ according
to (\ref{decomp}), then
$x_+$ cannot be fixed by $\phi$ in general.

%***************************************

%************************************
\section{Constant curvature immersions into pseudo-Riemannian space forms}
\label{part2}
In this section,  we first sketch the loop group construction of Ferus 
and Pedit \cite{feruspedit1996} for isometric immersions of 
space forms. We then evaluate these loop group maps along other
ranges of the spectral parameter, to obtain several equivalences
between isometric immersions into space forms and pseudo-Riemannian
space forms.

\subsection{The loop group formulation for isometric immersions
of space forms} \label{sketch}
Here is a brief outline of the formulation from \cite{feruspedit1996}.
In that work it was also shown how to construct many examples of
these maps using the AKS theory.
Let $M$ be a simply connected manifold of dimension $m$.
We first consider the case that the target space is a sphere.
The loop group maps are  elements of 
\beq \label{spaceform}
\frac{\hh_{{\rho} \mu}}{\hh_{{\rho} \mu}^0}  (M)_{-1}^1,
\eeq
where
$\hh = \Lambda G_\sigma$, 
 $G = SO(m+k+1,{\bf C})$, 
 $\sigma$ is
given by the equation (\ref{sigma}) for the involution $\sigma^0$ defining
a symmetric space $SO(m+k+1)/(SO(m) \times SO(k+1))$,
 namely 
\bdm
\sigma^0 := \textup{Ad}_P,
\edm
 for $P = \textup{diag}(I_{m}, -I_{k+1})$,
and $\mu$ is the negative  involution 
\beq  \label{mu}
(\mu X)(\lambda) := \textup{Ad}_Q(X(1/\lambda)),
\eeq
for $Q = \textup{diag}(I_{m +1}, -I_k)$.
Here $I_j$ is the $j\times j$ identity matrix.
Finally,  $\rho$ is one of three
reality conditions, described below,
 and $\hh_{\rho \mu}$ is the subgroup of $\hh$
fixed by both involutions $\rho$ and $\mu$.

There are essentially three cases for the induced Gauss curvature on
the immersion, which is constant for a fixed value of $\lambda$,
but varies with $\lambda$. 
These correspond to three different choices for the
reality condition ${\rho}$, and are displayed in Table~\ref{rangetable}.

\begin{table}[here]  
  \begin{tabular}{|c|c|c|}  \hline    
Reality condition & Parameter range & Induced Gauss curvature  \\ 

 $(\rho_1 X)(\lambda) :=  \overline{X(-\bar{\lambda})}$ & 
$\lambda \in i\real^*$ &
$c_\lambda \in (-\infty,0)$ \\

$(\rho_2 X)(\lambda) := 
  \overline{X(\bar{\lambda})}$ & $\lambda \in \real^*$ &
$c_\lambda \in (0,1]$ \\

$(\rho_3 X)(\lambda):=  \overline{X(1/\bar{\lambda})}$ & 
$\lambda \in S^1$ & 
$c_\lambda \in [1,\infty)$  \\ \hline 
  \end{tabular} 
\caption{Cases 1-3 for immersions into a sphere}
 \label{rangetable}
\end{table}

If $U$ is a subset of $M$, an \emph{adapted frame} for an 
immersion $f: U \to S^{m+k}$ is
defined to be a map 
$F = [e_1,...,e_m, f, \xi_1,....,\xi_k]: U \to SO(m+k+1)$,
 whose first $m$ and
last $k$ columns span the tangent and normal bundles to the image
$f(U)$.
The involutions $\sigma$ and $\mu$ mean that 
if $F$ is any representative of an element of
$\frac{\hh_{{\rho} \mu}}{\hh_{{\rho} \mu}^0}  (M)_{-1}^1$,
then, for a fixed value of $\lambda$, $F$ has the interpretation
as an  adapted frame for a map $f: M \to S^{m+k}$, provided the 
$(m+1)$'th column, $f$, is an immersion. The map $f$ is independent
of the choice of representative $F$,  because $\hh_\mu^0$ 
is just the subgroup $SO(m) \times {I} \times SO(k)$, which acts
on the right by fixing $f$ and changing the orthonormal frames for
the tangent and normal spaces. The involution $\mu$ also ensures
that the derivative $\dd f$ has no component in the directions
of any of $\xi_1,...,\xi_k$.

Notationally, we will not normally distinguish such a representative
$F$ from its equivalence class, and we will call either an
\emph{extended frame} for the family of immersions $f^\lambda$.
An element $F$ of $\hh_{{\rho} {\mu}}(M)_{-1}^1$ has a 
Maurer-Cartan form which looks like:
\beq \label{typeB}
F^{-1} \dd F  = \bbar {cc}
         \omega & \bbar {cc} (\lambda + \lambda^{-1}) \theta ~&~~ (\lambda - \lambda^{-1}) \beta \ebar \\
         -\bbar {cc} \epsilon (\lambda + \lambda^{-1}) \theta ~&~~ (\lambda - \lambda^{-1}) \beta \ebar^t  & \eta 
              \ebar,
\eeq
where
  the first row and column of $\eta$ are zero, and, 
in the spherical case, $\epsilon = 1$.
The 1-forms $\omega$ and $\eta$ are the connections of the tangent
and normal bundles respectively, $(\lambda + \lambda^{-1})\theta$ is 
the dual frame to our
tangent frame, and $(\lambda - \lambda^{-1})\beta$ is 
the second fundamental form.

A straightforward computation shows that $F^{-1} \dd F$
satisfying the Maurer-Cartan equation for all $\lambda$ is equivalent to
the integrability of $F$ at a single value of $\lambda$ plus the
extra conditions
\beqa
\dd \omega + \omega \wedge \omega = 4 \theta \wedge \theta^t \label{curvature},\\
\dd \eta + \eta \wedge \eta = 0.
\eeqa
The second equation says that the normal bundle is flat,
and equation (\ref{curvature}) says that the induced sectional 
curvature on the image of $f$ is 
\bdm
c_\lambda = \frac{4}{(\lambda + \lambda^{-1})^2},
\edm
which follows from the fact that the coframe 
 is $(\lambda + \lambda^{-1})\theta$. One then checks that as
$\lambda$ varies over the ranges $i \real^*$,  $\real^*$, and $S^1$,
 $c_\lambda$ varies 
over the intervals in Table~\ref{rangetable}.
The loop group map $F$ is well defined on $\cc^*$, but the map
$f$ is not an immersion at $\lambda = \pm i$, since the coframe
necessarily vanishes there.

If we allow our immersions to be degenerate at some points,
meaning the derivative drops rank, then any element of 
$\frac{\hh_{{\rho} \mu}}{\hh_{{\rho} \mu}^0}  (M)_{-1}^1$,
 corresponds to a family of isometric immersions.

Conversely, if $M$ is simply connected, then, once
 we fix the base point $p \in M$ at which our elements of
$\hh_{{\rho} {\mu}}(M)$ are normalised,  there is a unique element
of $\frac{\hh_{{\rho} \mu}}{\hh_{{\rho} \mu}^0}  (M)_{-1}^1$ associated
to a given isometric immersion with constant curvature in the
appropriate range of the table, with the exception of the 
limiting value $c_\lambda = 1$ in Cases 2 and 3.  This corresponds
to the values $\lambda = \pm 1$, at which point the second fundamental
form, $(\lambda - \lambda^{-1})\beta$, vanishes, so the immersion
is totally geodesic.  However, given a totally geodesic immersion,
we cannot insert $\lambda$ into its Maurer-Cartan form to obtain
the family, as we do not know what $\beta$ should be.

This converse statement was shown locally in \cite{feruspedit1996},
by choosing an adapted frame for the immersion.  A single
global adapted frame may not exist, but the equivalence class
in the  space 
$\frac{\hh_{{\rho} \mu}}{\hh_{{\rho} \mu}^0}  (M)_{-1}^1$
is nevertheless well defined:
%********************************
\begin{lemma} \label{welldefined}
Let $M$ be a simply connected manifold of dimension $m$,
and $f: M \to S^{m+k}$ an immersion
with flat normal bundle and induced constant curvature $c$,
with $0 \neq c \neq 1$, 
with the normalisation $f(p) = [0,...,1,0,...,0]$, the standard unit
vector $E_{m+1}$.  Then
there is a unique element $F \in 
\frac{\hh_{{\rho} \mu}}{\hh_{{\rho} \mu}^0}  (M)_{-1}^1$, where $\rho$
is the appropriate reality condition,
whose $(m+1)$'th column, evaluated at 
$\lambda_0 = \frac{1}{\sqrt{c}}(1+ \sqrt{1-c})$, is $f$.
\end{lemma}
%**************
\begin{proof}
For any point $q$ of $M$, there is a simply connected neighbourhood, $U_q$,
of $q$, which contains $p$, and an
 adapted frame, $F_q = [e_1,...,e_m,f,\xi_1,...,\xi_k]$, on $U_q$, 
normalised to the identity at $p$.  This can be obtained by parallel
translating the identity matrix along some path from $p$ to $q$, and
then extending to a neighbourhood of that path.
The Maurer-Cartan form of $F_q$ is 
\bdm
A_q := F_q^{-1} \dd F_q = \bbar 
{cc}
         \omega & \bbar {cc}  \theta ~&~~  \beta \ebar \\
         -\bbar {cc} \theta ~&~~  \beta \ebar^t  & \eta 
              \ebar.
\edm

Following \cite{feruspedit1996}, to insert the parameter
 $\lambda$, one multiplies $\theta$ by
$\frac{\sqrt{c}}{2}(\lambda + \lambda^{-1})$ and 
$\beta$ by $\frac{\sqrt{c}}{2\sqrt{1-c}}(\lambda - \lambda^{-1})$,
and then integrates on $U_q$,
with the initial condition $F(p) = I$, to get a representative
$F_q^\lambda \in \hh_{{\rho} \mu} (U_q )_{-1}^1$.

We need to check that for any other point $r$, $F_q^\lambda$ and
 $F_r^\lambda$ 
differ only by post-multiplication by maps into $\hh_{\rho \mu}^0$, on the
intersections of their domains of definition. 
Now at $x \in U_q \cap U_r$,
we have $F_r(x) = F_q(x) G(x)$ for some $G$ 
 which is smooth and takes values in 
$\hh_{\rho \mu}^0 = SO(m) \times 1 \times SO(k)$,
because $F_q$ and $F_r$ are both adapted frames. The matrix $G$
is of the form $\textup{diag}(A,1,B)$, and hence
\bdm
F_r^{-1} \dd F_r = \bbar 
{cc}
         A^t \omega A + A^t \dd A & \bbar {cc}  A^t \theta ~&~~ A^t \beta  B \ebar \\
         -\bbar {cc} A^t \theta ~&~~  A^t \beta B \ebar^t  & B^t \eta B + B^t \dd B
              \ebar.
\edm
It follows from the construction of the extended frame
 $F_r^\lambda$, that it has
the same Maurer-Cartan form at $x$ as $F_q^\lambda G$.
Since both functions are equal at $x$, and their Maurer-Cartan forms
agree, it follows that $F_r^\lambda$ and 
 $F_q^\lambda G$ agree wherever they are defined.
\end{proof}

The case where the target is hyperbolic space has the
same formulation, replacing the group $SO(m+k+1,\cc)$ with
$SO(m+k,1,\cc)$, defined here  to be the subgroup of
$GL(m+k+1,\cc)$ consisting of matrices 
which preserve the bilinear form given by
the matrix 
\bdm
J := \textup{diag}(I_m,-1,I_k).
\edm
In this case, $\epsilon = -1$ in (\ref{typeB}), and the corresponding 
induced curvatures in Table~\ref{rangetable} are all negated.

We will only need one of the cases for the hyperbolic space and
we therefore add to Cases 1-3 above, the following: \\
\noindent \textbf{Case 4}: 
$\Lambda SO(m+k,1,\cc)_{\sigma \mu \rho_2}(M)_{-1}^1$,
where $\sigma$, $\mu$ and $\rho_2$ are as before. Elements of
$\Lambda SO(m+k,1,\cc)_{\sigma \mu \rho_2}(M)_{-1}^1$ are isometric
immersions with flat normal bundle 
$M \to H^{m+k+1}$ with induced constant curvature $c_\lambda$ 
in the interval $[-1, 0)$.

%********************************
\subsection{Interpretation for other ranges of the spectral parameter}
The goal of this subsection is to identify for each of Cases 1-4 above,
 the different maps obtained for 
 values of the spectral parameter in all three ranges  $\real$, $i\real$ and
$S^1$.  This was partially investigated in \cite{branderdorf},
where the last row of Table~\ref{table1} below and the first row of
 Table~\ref{table2}
 were found.

Let $\real ^n_s$ denote the \emph{pseudo-Euclidean} space
$\real^n$ equipt with a metric $<,>$ with signature $(s,n-s)$, that is,
it is isometric to a space with metric $J= \textup{diag}(-I_s, I_{n-s})$.
Define the pseudo-Riemannian sphere, $S^n_s$, and pseudo-Riemannian
hyperbolic space, $H^n_s$, by
\beqas
S^n_s := \{ x \in \real_s^{n+1} ~:~ <x,x> = 1\},\\
H^n_s := \{ x \in \real_{s+1}^{n+1} ~:~ <x,x> = -1\}.
\eeqas
These two spaces are complete pseudo-Riemannian manifolds, both with signature
$(s,n-s)$, and with constant sectional curvatures 1 and $-1$ respectively
(see, for example, \cite{wolf}).

\begin{theorem} \label{corprop}
Let $G$ and $H$ denote the groups $SO(m+k+1,\cc)$ and
$SO(m+k,1,\cc)$ respectively. Let $F$ be an element of either
$\frac{\Lambda G_{\sigma \mu \rho_i}}{\Lambda G_{\sigma \mu \rho_i}^0}(M)_{-1}^1$, $i =1,...,3$, or
$\frac{\Lambda H_{\sigma \mu \rho_2}}{\Lambda H_{\sigma \mu \rho_2}^0}(M)_{-1}^1$, as described for Cases 1-4 above.  If the $(m+1)$'th column, $f$, of $F$ is
an immersion for all $\lambda \in \cc^* \setminus \{\pm i \}$, then, by evaluating $f$
for values of $\lambda$ in the ranges indicated in  
Tables~\ref{table1}-\ref{table4} below,
isometric immersions with flat normal bundle into the pseudo-Riemannian
space forms displayed are obtained.
In all cases, the induced metric on $M$  is positive
definite and with constant sectional curvature varying through  the range
indicated. 

 Conversely, if $M$ is simply connected, then
 any isometric immersion with flat 
normal bundle into one of the target spaces displayed,  with
constant sectional curvature in the corresponding range, apart
from the cases $c =1$ when the target space has positive curvature,
and $c=-1$ when the target space has negative curvature, belongs
to one of these families.

\begin{table}[here] 
  \begin{tabular}{|c|c|c|}  \hline
Parameter range & Induced Gauss curvature  & Target space \\ \hline
$\lambda \in i \real^* \setminus \{\pm i \}$  & $c_\lambda \in (-\infty,0)$ & $S^{m+k}$  \\
$\lambda \in  \real^*$  & $c_\lambda \in [-1,0)$ & $H^{m+k}_{k}$  \\
$\lambda \in S^1 \setminus \{\pm i\}$  & $c_\lambda \in (-\infty,-1]$ & $H^{m+k}$ \\
\hline
  \end{tabular}
\caption{Case 1: $~~$
$\frac{\Lambda G_{\sigma \mu \rho_1}}{G_{\sigma \mu {\rho_1}}^0}(M)_{-1}^1$}
\label{table1}
\end{table}

\begin{table}[here]
  \begin{tabular}{|c|c|c|}  \hline
Parameter range & Induced Gauss curvature  & Target space \\ \hline
$\lambda \in i \real^*\setminus \{\pm i\}$  & $c_\lambda \in (0,\infty)$ & $H^{m+k}_{k}$ \\
$\lambda \in  \real^*$  & $c_\lambda \in (0,1]$ & $S^{m+k}$ \\
$\lambda \in S^1\setminus \{\pm i\}$  & $c_\lambda \in [1,\infty)$ & $S^{m+k}_k$ \\
\hline
  \end{tabular}
\caption{Case 2: $~~$ 
$\frac{\Lambda G_{\sigma \mu \rho_2}}{G_{\sigma \mu {\rho_2}}^0}(M)_{-1}^1$}
 \label{table2} \end{table}

\begin{table}[here] 
  \begin{tabular}{|c|c|c|}  \hline
Parameter range & Induced Gauss curvature  & Target space \\ \hline
$\lambda \in i \real^*\setminus \{\pm i\}$  & $c_\lambda \in (0,\infty)$ & $H^{m+k}$ \\
$\lambda \in  \real^*$  & $c_\lambda \in (0,1]$ & $S^{m+k}_k$ \\
$\lambda \in S^1\setminus \{\pm i\}$  & $c_\lambda \in [1,\infty)$ & $S^{m+k}$ \\
\hline
  \end{tabular}
\caption{Case 3: $~~$
$\frac{\Lambda G_{\sigma \mu \rho_3}}{G_{\sigma \mu {\rho_3}}^0}(M)_{-1}^1$}
\label{table3}
\end{table}

\begin{table}[here]  
  \begin{tabular}{|c|c|c|}  \hline   
Parameter range & Induced Gauss curvature  & Target space \\ \hline
$\lambda \in i \real^* \setminus\{\pm i\}$  & $c_\lambda \in (-\infty,0)$ & $S^{m+k}_k$ \\
$\lambda \in  \real^*$  & $c_\lambda \in [-1,0)$ & $H^{m+k}$ \\
$\lambda \in S^1\setminus \{ \pm i\}$  & $c_\lambda \in (-\infty,-1]$ & $H^{m+k}_{k}$ \\
\hline
  \end{tabular}
\caption{Case 4: $~~$
$\frac{\Lambda H_{\sigma \mu \rho_2}}{H_{\sigma \mu {\rho_2}}^0}(M)_{-1}^1$}
\label{table4}  \end{table}
\end{theorem}
%___________________
\begin{proof} 
\noindent \textbf{Case 1} \\
Let $F \in \frac{\Lambda G_{\sigma \mu \rho_1}}{G_{\sigma \mu {\rho_1}}^0}(M)_{-1}^1$.  The first row of Table~\ref{table1} we already know from \cite{feruspedit1996}.
To get the second row, we need a reality condition along $\real$. Consider
$\phi := Ad_T : GL(m+k+1,\cc) \to GL(m+k+1,\cc)$, where
 $T:= \textup{diag}(iI_m,1,I_k)$. Now $\phi$ is an isomorphism between
$SO(m+k+1,\cc)$ and $SO(m,k+1,\cc)$, where the latter group is defined 
here to be the set of matrices preserving the bilinear form 
\bdm
\hat{J} := \textup{diag}(I_n, -I_{k+1}).
\edm
To verify this one checks that $A^t A =I$ is
 equivalent to $(\phi A)^t \hat{J} \phi  = \hat{J}$.

It is also easy to see that $\phi$ is a bijection between the two sets
$\frac{\Lambda G_{\sigma \mu \rho_1}}{G_{\sigma \mu {\rho_1}}^0}(M)_{-1}^1$
and $\frac{\Lambda SO(m,k+1,\cc)_{\sigma \mu \rho_2}}{SO(m,k+1,\cc)_{\sigma \mu {\rho_2}}^0}(M)_{-1}^1$, since $\phi$ commutes with both $\sigma$ and $\mu$ 
and if $F$ satisfies $\rho_1$ then
\beqas
\rho_2 (\phi F)(\lambda) &=& \overline{TF(\bar{\lambda})T^{-1}} \\
&=& \bar{T} P \overline{F(-\bar{\lambda})} P \bar{T}^{-1} \\
&=& \bar{T} P F(\lambda) P \bar{T}^{-1} \\
&=& (-T) F (\lambda) (-T)^{-1}\\
&=& \phi F(\lambda).
\eeqas
To get to the second line we used that fact that 
$\sigma F (\lambda) = \textup{Ad}_P F(-\lambda) = F(\lambda)$.

We now  interpret $\hat{F} := \phi(F)$.
The analysis of an element $\hat{F} \in
\frac{\Lambda SO(m,k+1,\cc)_{\sigma \mu \rho_2}}{SO(m,k+1,\cc)_{\sigma \mu {\rho_2}}^0}(M)_{-1}^1$ is similar to that of the case where the group is
$SO(m+k+1,\cc)$, explained above in Section \ref{sketch}: the Maurer-Cartan form has the expression
\bdm
\hat{F}^{-1} \dd \hat{F}  = \bbar {cc}
         \omega & \bbar {cc} (\lambda + \lambda^{-1}) \theta ~&~~ (\lambda - \lambda^{-1}) \beta \ebar \\
         \bbar {cc} (\lambda + \lambda^{-1}) \theta ~&~~ (\lambda - \lambda^{-1}) \beta \ebar^t  & \eta 
              \ebar,
\edm
where the first row and column of $\eta$ are zeros.
This differs from the expression (\ref{typeB}) only by a minus
sign in the lower left corner.
$\hat{F}$ is real for values of $\lambda$ in $\real$.
For such values of $\lambda$,
 we take the $(m+1)$'th column of $\hat{F}$ as our map
$f$, then, by the definition of $SO(m,k+1)$, we have $f^t \hat{J} f = -1$, so
$f$ takes values in $H^{m+k}_{k+1}$.  The zeros in the first row and
column of $\eta$ say that 
the tangent space to the image of $f$ lies in the span of the first $m$
columns, so that if $f$ is an immersion then the first $m$
columns, $\{e_1,...,e_m\}$, are a frame for the tangent bundle. 
Since $\hat{F}^t \hat{J} \hat{F} = \hat{J}$, it follows from the form
of $\hat{J}$ that $e_i^t \hat{J} e_j = \delta_{ij}$, so 
the induced metric is positive definite. 
As before, the integrability condition implies the equations,
\beqa
\dd \omega + \omega \wedge \omega = -4 \theta \wedge \theta^t,\\
\dd \eta + \eta \wedge \eta = 0.
\eeqa
which imply that the normal bundle is flat and the induced Gauss curvature is:
\bdm
c_\lambda = \frac{-4}{(\lambda + \lambda^{-1})^2},
\edm
which varies over the interval $[-1,0)$ as $\lambda$ varies over $\real^*$.
The converse argument is also similar to the spherical case of 
Lemma \ref{welldefined}.

The third row of Table~\ref{table1} was obtained in \cite{branderdorf}, in a 
similar fashion, using $T = \textup{diag}(iI_m, 1, iI_k)$ and
$\hat{J} = \textup{diag}(I_m,-1,I_k)$.

\noindent \textbf{Case 2} \\
The second row is given in \cite{feruspedit1996}. To get the first row,
 proceed as in Case 1, using $T = \textup{diag}(iI_m,1,I_k)$ and 
$\hat{J} = \textup{diag}(I_m,-1,-I_k)$. For the third row, use 
$T = \textup{diag}(iI_m,i,I_k)$ and 
$\hat{J} = \textup{diag}(I_m,1,-I_k)$.\\

\noindent \textbf{Case 3} \\
The third row is given in \cite{feruspedit1996}. To get the first row,
again proceed as in Case 1, using $T = \textup{diag}(iI_m,1,iI_k)$ and 
$\hat{J} = \textup{diag}(I_m,-1,I_k)$. For the second row, use 
$T = \textup{diag}(iI_m,i,I_k)$ and 
$\hat{J} = \textup{diag}(I_m,1,-I_k)$.\\

\noindent \textbf{Case 4} \\
The second row we know from \cite{feruspedit1996}. To get the first row,
 use  $T = \textup{diag}(I_m,i,iI_k)$ and 
$\hat{J} = \textup{diag}(I_m,1,-I_k)$. For the third row, use 
$T = \textup{diag}(I_m,1,iI_k)$ and 
$\hat{J} = \textup{diag}(I_m,-1,-I_k)$.\\
\end{proof}

\begin{remark} \label{other}
It is clear from the preceding proof that isometric immersions 
with flat normal bundle of a constant curvature Riemannian manifold 
$M^m$ into either $S_l^{m+k}$ or $H_l^{m+k}$, for any $0 \leq l \leq k$
can be treated similarly, by starting with the group which preserves
the bilinear form $J = \textup{diag}[I_m, \epsilon ,-I_l,I_{k-l}]$,
where $\epsilon = \pm 1$. 
\end{remark}
%************** 
\begin{example} \label{example}
Here is a simple  example from Case 3. Consider the family of maps 
$F_\lambda : \real ^2 \to G = SO(4,\cc)$ which takes $(u,v) \in \real^2$
to the matrix
\bdm
 \bbar {cccc}
     \cos(u) & -\sin(u) \sin(v) & a \sin(u) \cos(v) & b \sin(u) \cos(v) \\
      0  & \cos(v)               & a \sin(v)         & b \sin(v)  \\
 -a \sin(u) & -a\cos(u)\sin(v) & a^2 \cos(u)\cos(v) + b^2 & 
               ab(\cos(u)\cos(v)-1)\\
     -b\sin(u) &  -b\cos(u)\sin(v) & ab(\cos(u)\cos(v)-1) &
                 b^2\cos(u)\cos(v)+a^2   \ebar,
\edm
where 
\bdm
a= \frac{1}{2}(\lambda + \lambda^{-1}),  \hspace{1cm} 
b= \frac{i}{2}(\lambda - \lambda^{-1}).
\edm

The Maurer-Cartan form of $F_\lambda$ is
\bdm
F_\lambda ^{-1}\dd F_\lambda = \bbar {cccc}
   0 & -\sin(v) \dd u & a \cos(v) \, \dd u & b \cos(v) \, \dd u \\
  \sin(v) \dd u & 0 & a \, \dd v &  \, b \dd v \\
  - a \cos(v) \, \dd u & -a \, \dd v & 0 & 0 \\
  -b \cos(v) \,  \dd u & -b \, \dd v & 0 & 0
  \ebar.       
\edm
Now $F_\lambda^{-1}\dd F_\lambda$
 is fixed by $\sigma$, $\mu$ and $\rho_3$, so it
takes values in the Lie algebra of $\Lambda G_{\sigma \mu \rho_3}$,
and $F_\lambda (0,0) = I \in \Lambda G_{\sigma \mu \rho_3}$. 
Therefore $F_\lambda$ is a map into $\Lambda G_{\sigma \mu \rho_3}$,
and, since its Maurer-Cartan form has top and bottom degree 1 and -1
respectively, it represents an element of 
$\frac{\Lambda G_{\sigma \mu \rho_3}}{G_{\sigma \mu \rho_3}^0}(\real^2)_{-1}^1$. 
Thus, according to Table~\ref{table3}, if the third column of $F_\lambda$, namely
\bdm
f^\lambda (u,v) = \bbar {c} \frac{1}{2}(\lambda + \lambda^{-1}) \sin (u) \cos(v) \\
         \frac{1}{2} (\lambda + \lambda^{-1}) \sin(v) \\
        \frac{1}{4}(\lambda + \lambda^{-1})^2 \cos(u) \cos(v) - \frac{1}{4}((\lambda - \lambda^{-1})^2 \\ 
         \frac{i}{4}((\lambda + \lambda^{-1})(\lambda - \lambda^{-1})(\cos(u)\cos(v) -1) \ebar =: \bbar {c} f_1\\ f_2 \\ f_3\\ f_4 \ebar,
\edm 
 is an immersion, then, for a value of $\lambda$ in $S^1$, it is an immersion
into $S^3$ with constant Gauss curvature greater or equal to 1.
The dual frame for $f^\lambda$ is given by 
\beq \label{theta}
\theta = \bbar {c} \frac{1}{2}(\lambda + \lambda^{-1}) \cos(v) \, \dd u \\
           \frac{1}{2}(\lambda + \lambda^{-1}) \, \dd v \ebar,
\eeq
and so, if $\lambda \neq \pm i$, then $f^\lambda$ is
immersive away from the 
 degenerate coordinate lines $\cos(v) =0$.
In fact $f^\lambda$ is a deformation, through a family of isometrically
embedded spheres, of the totally geodesic embedding of
$S^2$ into $S^3$ given by 
\bdm
f(u,v) = \bbar {cccc}  \sin (u) \cos(v), &
          \sin(v), &
        \cos(u) \cos(v), &
         0 \ebar^t,
\edm
 which is achieved at $\lambda = 1$.

To obtain the isometric immersions from the first two lines of
 Table~\ref{table3},
we need to apply the transformations $Ad_T$ given in the proof first:
using $T = \textup{diag}(iI_2, 1, i)$,
 the third column of 
$\textup{Ad}_T(F_\lambda)$ is
\bdm
\hat{f}^\lambda (u,v) = [if_1, if_2, f_3, if_4]^t,
\edm
and this is indeed real for values of $\lambda$ along $i\real$, and
gives a family of embeddings  of a  sphere into $H^3$, with constant sectional
curvature $c_\lambda = \frac{-4}{(\lambda + \lambda^{-1})^2}$, which
varies through the range $(0,\infty)$.

Finally, to get the immersion in the  second row of Table~\ref{table3}, 
we use the matrix $T = \textup{diag}(iI_2,i,1)$ to obtain the family 
\bdm
\tilde{f}^\lambda (u,v) = [f_1, f_2, f_3, -if_4]^t.
\edm
This is real for real values of $\lambda$. At $\lambda = 1$ it agrees with 
$f^\lambda$, being the same sphere embedded in the plane $\real^3$.
As $\lambda$ varies over $\real^*$, the immersion moves  
through a family of isometrically embedded spheres, with constant curvature
$c_\lambda \in (0,1)$, in the de Sitter space $S^3_1$.
\end{example}
\begin{remark}  Example \ref{example} is not typical. 
As the parameter $\lambda$ varies,  one should not
normally expect an embedding to remain an embedding, 
but, rather, only an immersion.
\end{remark}

%****************************************************************

\subsection{Relations between the immersions from different ranges
of the spectral parameter}
The deformation parameter  will not generally appear
in the maps $f^\lambda$ in such a simple manner as occurred in
Example \ref{example}.   However, at the level of the Maurer-Cartan
form, it is always the same, and therefore one has the following 
global result concerning the maps $f^\lambda$, $\hat{f}^\lambda$ and
$\tilde{f}^\lambda$ in general:

\begin{proposition} \label{metriclemma}
Let $f_j^{\lambda_j}$, $j=1,..,2$, 
be two maps from any one of  the Tables~\ref{table1}-\ref{table4}, 
obtained from the $(m+1)$'th column of $F$ or $\textup{Ad}_T F$, as
described in Theorem \ref{corprop},
evaluated at two given points
 $\lambda_j  \in (i\real^*  \cup
\real ^* \cup  S^1) \setminus \{\pm i \})$.
Then:
\begin{enumerate}
\item
If $f_1^{\lambda_1}$ is immersive at a 
point $x \in M$,  then so is $f_2^{\lambda_2}$.
\item
If  $f_1^{\lambda_1}$ is a complete immersion then
so is $f_2^{\lambda_2}$.
\end{enumerate}
\end{proposition}
\begin{proof}
The  matrix $T$ used 
 to go between $f_1$ and $f_2$ is among the following
list: $I$, $\textup{diag}(iI_m,1,I_k)$, $\textup{diag}(iI_m,1,iI_k)$,
$\textup{diag}(iI_m,i,I_k)$, $\textup{diag}(I_m,i,iI_k)$
and $\textup{diag}(I_m,1,iI_k)$, together with their compositions and
inverses.  It follows that the Maurer-Cartan forms of
$F$ and $\textup{Ad}_T F$ are both of the form:
\bdm
 \bbar {cc}
         \omega & \bbar {cc} (i)^r(\lambda + \lambda^{-1}) \theta ~&~~ (i)^s(\lambda - \lambda^{-1}) \beta \ebar \\
         \bbar {cc} \pm(i)^r(\lambda + \lambda^{-1}) \theta ~&~~ \pm(i)^s(\lambda - \lambda^{-1}) \beta \ebar^t  & \eta 
              \ebar,
\edm
for some fixed real matrix valued 1-forms 
$\omega$, $\theta$, $\beta$ and $\eta$. Thus the coframe for
 $f_j^{\lambda_j}$, obtained from the first $m$ components of
the $(m+1)$'th column, is just a non-zero constant $k_j$ times the 
fixed column vector valued 1-form $\theta$. The relevant reality
 condition ensures that
this constant is real.  The condition for $f_j^{\lambda_j}$ to be
an immersion is that the coframe consist of $m$ linearly independent 
1-forms, which proves the first part of the proposition.

For completeness, if $f_j^{\lambda_j}$ is an immersion with
 coframe 
 $k_j \theta = k_j[\theta_1,....,\theta_m]^t$,
then the induced metric is
\bdm
k_j^2(\theta_1^2 + ... + \theta_m^2).
\edm
Thus the induced metrics for the two immersions are positive constant
multiples of each other, and hence completeness is equivalent
for them.
\end{proof}

%--------------
\begin{corollary}  \label{apcor}
Within any one of Tables~\ref{table1}-\ref{table4} of Theorem \ref{corprop}, the
existence problem for an isometric immersion
with flat normal bundle $M_c^m \to N_{\tilde{c}}^{m+k}$, $c \neq \tilde{c}$,
 where $M_c^m$ is a
complete simply connected $m$-dimensional space form of constant
curvature $c$ in one of the appropriate intervals,
 and $N_{\tilde{c}}^{m+k}$
is the corresponding target space of constant curvature $\tilde{c}$, 
is equivalent throughout the table. 
\end{corollary}
\begin{proof} This follows from the converse part of Theorem
\ref{corprop}, together with Proposition \ref{metriclemma}.
\end{proof}

%************************
\subsection{Applications}
An interesting application of Corollary \ref{apcor} is to generalisations
of the well known theorem of Hilbert that $H^2$ cannot be globally immersed
into Euclidean space $E^3$ \cite{hilbert1}.
\begin{corollary} \label{apcor2}
Let $c$ be a negative real number. 
The problems of globally isometrically immersing the $m$-dimensional
simply connected space form $Q_c^m$ into the space forms $\tilde{Q}_{\tilde{c}}^{2m-1}$, for $c < \tilde{c}$, with $0 \neq \tilde{c}$, are all
equivalent.
\end{corollary}
\begin{proof}
For $c < \tilde{c}$, the normal bundle is automatically flat in this
codimension \cite{moore1}. Thus our problem is in the realm of Table~\ref{table1}
of Theorem \ref{corprop}. The case $\tilde{c} >0$ belongs to 
the first row of Table~\ref{table1}, after rescaling the sphere so that
$\tilde{c} =1$.
 The case $\tilde{c} < 0$ fits into the third line of
the table, after rescaling so that $\tilde{c} = -1$.
Hence,  Corollary \ref{apcor} implies the result.
\end{proof}
As mentioned in the introduction, Corollary \ref{apcor2} is a known
result, but our proof does not depend on special coordinates.

Another application is:
\begin{corollary}  \label{apcor3}
\begin{enumerate}
\item
There is no global isometric immersion with flat normal bundle 
of a
 sphere $S^m(R)$, of dimension $m$ and
any radius $R$, into the pseudo-Riemannian hyperbolic space
$H^{2m-1}_{m-1}$.
\item
There is no global isometric immersion with flat normal bundle
 of a
 sphere $S^m(R)$, of dimension $m$ and 
radius $R<1$, into the pseudo-Riemannian sphere $S_{m-1}^{2m-1}$.
\end{enumerate}
\end{corollary}
\begin{proof}
These follow from JD Moore's proof \cite{moore1} that a sphere
$S^m(R)$ of radius $R >1$, or equivalently  Gauss curvature less 
than 1, cannot be globally immersed into $S^{2m-1}$. 
Together with Corollary \ref{apcor}, this says
 that complete immersions are not possible in
Table~\ref{table2} of Theorem \ref{corprop}, which accounts
for both cases of this corollary.
\end{proof}
\begin{remark} In the special case that $m=2$, then the normal bundle
is flat, so Corollary \ref{apcor3} reproduces  the result
of Li \cite{li}, that there is no isometric immersion of a 2-sphere of 
constant curvature greater than 1 into the de Sitter space $S^3_1$.
\end{remark}

%*************************************************************************************

\section{The DPW method} \label{dpw}
The next goal of this paper is to relate the loop group formulation
of isometric immersions to that of pluriharmonic maps. For this we
will need the generalised DPW method.
The DPW method was first used in \cite{dorfmeisterpeditwu} to produce
 harmonic maps from a Riemann surface into a symmetric space from holomorphic data.
It was extended to pluriharmonic maps in \cite{dorfech}.
The main idea of the method was shown to be extendable to somewhat arbitrary connection 
order $(_a^b)$ maps in \cite{branderdorf}, to which we refer the
reader for more details of the following sketch.

Let $\mathcal{H}$ be a Birkhoff decomposable subgroup of $\Lambda G$,
and $a\leq 0 \leq b$ be extended integers. The generalised DPW method
 gives the following bijection:
\beqa \label{bcor}
F \in \frac{\mathcal{H}}{\mathcal{H}^0}(M)_a^b ~~ \leftrightarrow  &
 F_+ \in \mathcal{H}(M)_1^b, \\
 & F_- \in \mathcal{H}(M)_a^{-1},  \nonumber
 \eeqa
which holds provided elements to be factored are 
in the big cell of $\Lambda G$.
The maps $F_+$ and $F_-$ are simply the left factors in the left
and right 
Birkhoff decompositions 
\beq \label{bdec}
F =F_+G_- = F_-G_+. 
\eeq
The method was used in this form in \cite{todathesis} for pseudospherical
surfaces in $\real^3$, in which case the
functions $F_+$ and $F_-$ were functions of independent variables,
simplifying the problem.

If $\tau$ is a negative involution of $\mathcal{H}$,
and $F \in \frac{\mathcal{H}_\tau}{\mathcal{H}_\tau^0}(M)_a^b$,
then it follows that $a = -b$ and,  applying
$\tau$ to  (\ref{bdec}),  we deduce from uniqueness of the
Birkhoff factorisation that
we must have $F_- = \tau F_+$. In fact it requires some work to 
prove the $\leftarrow$ side of the correspondence here, but what 
one has is
\beq  \label{taucor}
 F \in \frac{\mathcal{H}_\tau}{\mathcal{H}_\tau^0}(M)_{-b}^b
 ~~ \leftrightarrow  ~~ F_+ \in \mathcal{H}(M)_1^b. 
  \eeq
For the main purpose of what follows, we only need to know that the
bijection $(\ref{taucor})$ always holds in a neighbourhood of the
identity, although it is true that if
 $\tau$ is the $S^1$ reality condition (\ref{tau}), where $\tau^0$
defines a compact real form of $G$, then the $\leftarrow$
correspondence of (\ref{taucor}) is global on $M$, as it is
constructed from an Iwasawa splitting of the loop group, which
holds globally.

%**********************************************************

\section{Curved flats} \label{curvedflats}

Curved flats were defined in \cite{feruspedit1996II} as follows:
let $U/K$ be a semisimple symmetric space defined by the 
commuting involutions
$\sigma^0$ and $\rho^0$ of a complex semisimple Lie group $G$,
where $U$ is the fixed point set of the reality condition 
$\rho^0$ and $K$ is the fixed point set of both involutions.
Let  $M$ be a connected  manifold of dimension $m$. The map
$f: M \to U/K$ is a curved flat if $f^*R = 0$ as a 2-form on
$M$, where $R$ is the curvature tensor of $U/K$. 

 Here we are 
principally interested in the loop group formulation which defines
a \emph{family} of curved flats $f^\lambda$, parameterised by 
$\lambda$  in the nonzero real numbers $\real ^*$.
To define these, extend $\rho^0$ to
an involution $\rho$ of $\Lambda G$ by the formula (\ref{rho}),
so that elements of $\Lambda G_\rho$ 
are in $U$ for $\lambda \in \real^*$.
We also extend $\sigma^0$ to an involution of $\Lambda G$ by the
formula (\ref{sigma}).

Let $\hh := \Lambda G_\sigma$, the fixed point subgroup of $\sigma$.
Note that $\rho$ and $\sigma$ commute, and we define
$\hh_\rho := \Lambda G_{\rho \sigma}$ to be the subgroup of elements of
 $\Lambda G$ fixed
by both involutions.  It is shown in \cite{feruspedit1996II} that
a (family of lifts into U of) curved flats is just an element of 
$\hh _\rho (M)_0^1$, that is, a map $F$ from $M$ into
$\hh _\rho$ whose Maurer-Cartan form has the
expansion
\bdm
A^\lambda := F_\lambda^{-1} \dd F_\lambda = A_0 + A_1 \lambda.
\edm
The involution $\sigma$ enforces that $A_0$ and $A_1$ are in the
$+1$ and $-1$ eigenspaces, $\mathfrak{k}$ and $\mathfrak{p}$,
 respectively of $\sigma^0$, and the 
\emph{curved flat equations} for $F$, namely 
\beqa
\dd A_0 + A_0 \wedge A_0 = 0, \label{a0int} \\
\dd A_1 + A_0 \wedge A_1 + A_1 \wedge A_0 = 0, \nonumber \\
A_1 \wedge A_1 = 0, \nonumber
\eeqa
are equivalent to the fact that $A^\lambda$ satisfies the 
Maurer-Cartan equation $\dd A + A \wedge A = 0$ for all $\lambda$.
This is  the integrability condition for $F$, so it must hold.

In fact, since $A_0$ is in $\mathfrak{k}$, the Lie algebra of $K$,
 and is itself integrable by (\ref{a0int}), 
we can gauge away this term by right multiplication by  a map into $K$.
In other words, the same family of maps into $U/K$ is represented
by a map $\hat{F}_\lambda$ whose Maurer-Cartan form has the expansion 
\bdm
\hat{A}^\lambda = \hat{A}_1 \lambda.
\edm
We therefore make the following definition:
\begin{definition} \label{curvedflatdef}
Let $G$, $U$, $K$, $\sigma$ and $\rho$ be as above. A 
\emph{(normalised) extended curved flat} from $M$ into $U/K$
 is an element of $\hh _\rho (M)_1^1$.
\end{definition}

%**********************************************
\section{Pluriharmonic maps into symmetric spaces} \label{pluri}
Harmonic maps were first studied in the loop group setting by
Uhlenbeck \cite{uhlenbeck1989}, and this was extended to 
pluriharmonic maps by Ohnita and Valli \cite {ohnitavalli}. For more
details on the formulation described here, as well as a discussion of
methods to produce \emph{finite type} examples, the reader could
consult \cite{bfpp}. The geometrical interpretation of the spectral
parameter deformation is described in \cite{dorfech}. We will proceed
directly to the loop group formulation here.

Let $\uhat / \khat$ be a semisimple symmetric space given 
by the involution $\sigma^0$ and reality condition $\hat{\tau}^0$
of $G$. In a later section we will assume that $G$ and
$\sigma^0$ are those given in Section \ref{curvedflats}, while
the reality conditions $\hat{\tau}^0$ and $\rho^0$ will be different.

We extend $\sigma^0$ to $\Lambda G$ again by the formula (\ref{sigma}),
but this time we extend our reality condition $\hat{\tau}^0$ in a 
different way, by the rule (\ref{tau}),
so that elements of $\Lambda G_{\hat{\tau}}$ are $\uhat$-valued 
 for unitary values of $\lambda$.

Let $\mcc$ be a simply connected $m$  dimensional complex manifold. As before, 
let $\hh := \Lambda G_\sigma$.
In \cite{bfpp} it is shown that an extended lift for a pluriharmonic
map from $\mcc$ into $\uhat / \khat$ is given by an element $F$ of
$\hh_{\hat{\tau}}(\mcc)_{-1}^1$, with one additional property,
namely, that if one expands the Maurer-Cartan form of $F$,
\beq \label{aa}
A^\lambda := F_\lambda^{-1} \dd F_\lambda =
 A_{-1}  \lambda^{-1} + A_0
 +  A_1 \lambda^{1},
\eeq
and
\bdm
A_i = A_i^\prime + A_i^{\prime \prime}
\edm
is the decomposition of the 1-form $A_i$ into its (1,0) and
(0,1) components  with 
respect to a complex basis for the cotangent space,
 then 
\bdm
A_1^{\prime \prime}  = 0.
\edm
Together with the fact that $A^\lambda$ is fixed by ${\hat{\tau}}$ this
also means that $A_{-1}^{\prime} = 0$ 
and $A^{\prime \prime}_{-1} = \overline{A^{\prime}_1}$, so that we have 
an expression
\bdm
A^\lambda := \overline{A^\prime_{1}}  \lambda^{-1}
    + A_0 + A^\prime_{1} \lambda.
\edm
The extended pluriharmonic map into $\uhat / \khat$ associated to an extended 
lift $F \in \hh_{\hat{\tau}}(\mcc)_{-1}^1$ is its equivalence class
modulo right multiplication by a map into $\khat$, which is just the subgroup
$\hh_ {\hat{\tau}}^0$ of constant loops, hence we make the following definition:

\begin{definition} \label{pluridef}
Let $G$, $\sigma$, ${\hat{\tau}}$, $\uhat$ and $\khat$ be as above.
An \emph{extended  pluriharmonic map} 
from $\mcc$ into $\uhat / \khat$
is an element $F$ of $\frac{\hh _{\hat{\tau}}}{\hh_{\hat{\tau}}^0} (\mcc)_{-1}^1$, 
with the property that if
$F_\lambda^{-1} \dd F_\lambda$ has the expansion (\ref{aa}),
then 
\beq \label{pluricond}
A_1^{\prime \prime} = 0.
\eeq
We will denote the set of these by 
$\mathcal{P}\frac{\hh _{\hat{\tau}}}{\hh_{\hat{\tau}}^0} (\mcc)_{-1}^1$
\end{definition}
\begin{remark}  \label{pluriremark}
 For every element 
$F \in \mathcal{P}\frac{\hh _{\hat{\tau}}}{\hh_{\hat{\tau}}^0} (\mcc)_{-1}^1$
there is a unique pluriharmonic map $f: \mcc \to \uhat/\khat$ obtained by
evaluating $F$ at $\lambda = 1$, and vice versa \cite{dorfech}.
\end{remark}

%*****************************
\subsection{DPW for pluriharmonic maps} \label{subsdpw}
Now if one applies the splitting described in Section \ref{dpw}
 to an extended
 pluriharmonic map as given in Definition \ref{pluridef} then it
is straightforward to check that the condition (\ref{pluricond})
implies that the map $F_+$ on the right hand side of (\ref{taucor})
is \emph{holomorphic} in the $\mcc$ variables,
 that is $\overline{\partial}F_+ = 0$.  Conversely,
one can show that if one starts with an element 
$F_+ \in \hh (\mcc)_1^1$ which is holomorphic on $\mcc$,
then the map $F$ given by the left hand side of (\ref{taucor}) is
pluriharmonic.   Now even though the DPW correspondence (\ref{taucor}) only
holds on the big cell in general, we can always renormalise our extended pluriharmonic
map $F$ at any point $q$ by premultiplying it by $F^{-1}(q,\lambda)$. 
This is again a pluriharmonic map,  and in the big cell on a neighbourhood of $q$,
and therefore we have an  alternative global characterisation
of pluriharmonic maps.  We  summarise this as:
\begin{proposition} \label{reform}
Let $G$, $\sigma$, ${\hat{\tau}}$, $\uhat$ and $\khat$ be given as in
Definition \ref{pluridef}.  Suppose $F$ is an element of
 $\frac{\hh_ {\hat{\tau}}}{\hh_{\hat{\tau}}^0}   (\mcc)_{-1}^1$.
Then $F$ is an extended pluriharmonic map if and 
only if the corresponding $F_+ \in \hh (\mcc)_1^1$ from
 the right hand side of  (\ref{taucor}) is holomorphic
on $\mcc$.
\end{proposition}
\begin{remark} Note that this test is understood to be applied locally by renormalizing
$F$.  The holomorphic function $F_+$ is not in general defined globally.
\end{remark}

%
%
%******************************************************************
%******************************************************************
%
\section{Pluriharmonic maps constructed from analytic curved flats} \label{pluriflats}
Let $G$, $U$, $K$, $\sigma$, $\rho$, $\uhat$, $\khat$ and ${\hat{\tau}}$
be as in Sections \ref{curvedflats} and \ref{pluri}, and assume
that the reality conditions $\rho^0$ and ${\hat{\tau}}^0$ commute, so
that the  extended involutions given by (\ref{rho}), (\ref{sigma}) and
(\ref{tau}) commute also.
Let $M$ be a connected paracompact real analytic manifold of dimension $m$.
By taking an atlas of $M$ and analytically extending the transition 
functions,  we can 
embed $M$ as a totally real submanifold of some complex manifold $M_\cc$
of complex dimension $m$.  The following theorem always holds at least
locally on $M$, and globally if $\khat$ is compact.

\begin{theorem} \label{theorem1}
Let $f_+: M \to U/K$ be a real analytic curved flat, represented by
the extended family 
\bdm
F_+ \in \hh_\rho (M)_1^1.
\edm
Then there exists an open submanifold $M_\epsilon$ of $M_\cc$, containing
$M$, and a unique pluriharmonic map
 $\hat{f}: M_\epsilon \to \uhat/\khat$,
represented by 
\bdm
\hat{F} \in \mathcal{P}\frac{\hh_{\hat{\tau}}}{\hh_{\hat{\tau}}^0}(M_\epsilon)_{-1}^1,
\edm
such that the restriction of $\hat{F}$ to $M$ satisfies the reality condition $\rho$.
More precisely we have the following 
correspondence from (\ref{taucor}):
\beq  \label{realpluri}
F = \hat{F}|_M \in \frac{\hh_{\rho {\hat{\tau}}}}{\hh_{\rho {\hat{\tau}}}^0}(M)_{-1}^1
 ~~ \leftrightarrow ~~ F_+ \in \hh_{\rho}(M)_1^1.
\eeq
\end{theorem}

\begin{proof}
We need only to show that there is a holomorphic extension  
$\hat{F}_+ \in \hh (M_\epsilon)_1^1$ of $F_+$, for
some open $M_\epsilon$ containing $M$.  Then the DPW correspondence
(\ref{taucor}) together with Proposition \ref{reform} 
gives us the required 
$\hat{F} \in  \mathcal{P}\frac{\hh_{\hat{\tau}}}{\hh_{\hat{\tau}}^0}(M_\epsilon)_{-1}^1$. 
 To see that $\hat{F}$ restricts on $M$ to an element of
$\frac{\hh_{\rho {\hat{\tau}}}}{\hh_{\rho {\hat{\tau}}}^0}(M)_{-1}^1$, observe
that $\hat{F}|_M$ is just the object obtained by applying the 
DPW correspondence (\ref{taucor}) to $F_+$ itself, since this is done
pointwise on $M_\epsilon$.  Since $\rho$ is a positive involution, $\hh_\rho$ 
is Birkhoff decomposable, and so  $\hat{F}|_M$ is 
fixed by $\rho$.
 The existence of $\hat{F}_+$ is shown in 
Proposition \ref{complexifyingprop} below.
\end{proof}

%***********************
\subsection{Complexifying real analytic curved flats}
To complete the proof of Theorem \ref{theorem1}, we need to show
there is a holomorphic extension of $F_+$. 
%******
\begin{proposition} \label{complexifyingprop}
Let $F_+$ be a real analytic element of $\hh_\rho(M)_1^1$, as above.
There exists an open submanifold $M_\epsilon$ of $\mcc$, containing
$M$, such that $F_+$ has a unique holomorphic extension to an element
$\hat{F}_+ \in \hh(M_\epsilon)_1^1$.
\end{proposition}
\begin{proof}
Consider the Maurer-Cartan form 
\bdm
A_+ := F_+^{-1} \dd F_+ = \eta \lambda,
\edm
where $\eta$ is an analytic $\mathfrak{p}$-valued 1-form on
$M$ satisfying the curved flat equations:
\beqa \label{cf1} 
\dd \eta = 0, \\
\eta \wedge \eta = 0, \label{cf2} 
\eeqa
which are equivalent to the Maurer-Cartan equation for $A_+$.
Let $\mathfrak{p}_\cc$ denote the complexification of  $\mathfrak{p}$.
%---------------lemma1------------
\begin{lemma}  \label{lemma1}
There exists an open submanifold,
$M_\epsilon$, of $\mcc$, containing 
$M$, such that $\eta$ has a unique analytic extension to
a $\mathfrak{p}_\cc$-valued holomorphic 1-form $\eta^\cc$ on  $M_\epsilon$,
satisfying the curved flat equations (\ref{cf1}) and
 (\ref{cf2}).
\end{lemma}

Let $z^j = x^j + iy^j$ be local coordinates on  $\mcc$. Then $\eta$
has the expression
\bdm
\eta = \sum_j \eta_j \dd x^j,
\edm
where each $\mathfrak{p}$ valued function $\eta_j$ 
is analytic in $x^1$,...,$x^m$
on $M$. 
 By standard theory of power series,
there is a  
neighbourhood $M _\epsilon$ of $M$ in $\mcc$ to which
each function $\eta _j$ has a holomorphic extension, 
$\hat{\eta}_j$, which takes values in $\mathfrak{p}_\cc$,
and this extension
is unique. Define $\hat{\eta} := \sum \hat{\eta}_k \dd z^k$.
We  assumed our manifold $\mcc$ was constructed from $M$ by
analytically extending the transition functions defining $M$,
from which it follows that $\hat{\eta}$ is a well defined
global holomorphic extension of $\eta$.

Let $d = \partial + \bar{\partial}$ be the usual decomposition
of  the $d$ operator into holomorphic and antiholomorphic parts.
Now  $\bar{\partial}_k \hat{\eta}_j = 0$, so we have
\bdm
\dd \hat{\eta} = \sum_{j<k} (\partial_j \hat{\eta} _k - 
  \partial_k \hat{\eta}_k)\dd z^j \wedge \dd z^k,
\edm
and the term $(\partial_j \hat{\eta}_k - \partial_k \hat{\eta}_k)$ is just the 
analytic extension of $(\frac{\partial \eta_k}{\partial x^j}
   - \frac{\partial \eta_j}{\partial x^k}$), which vanishes by (\ref{cf1}).
The argument for (\ref{cf2}) is analogous, and this proves the lemma.

Now denote by  $\hat{A}_+$ the family of 1-forms $\hat{\eta} \lambda$.
By Lemma \ref{lemma1}, $\hat{A}_+$ satisfies the Maurer-Cartan equations
for all $\lambda$, and therefore
integrates to a map $\hat{F}_+ \in \hh (\widetilde{M_\epsilon})_1^1$,
where  $\widetilde{M_\epsilon}$ is the universal cover of $M_\epsilon$.
$\hat{F}_+$ is uniquely
determined by the normalisation at some base point in $\pi^{-1}(p)$,
where $p$ is the base point of $M$, and  $\pi$ is the projection 
$\widetilde{M_\epsilon} \to M_\epsilon$.
Since  $\hat{\eta}$ is holomorphic on $\widetilde{M}_\epsilon$, 
so is $\hat{F}_+$.

Finally, we need to show that $\hat{F}_+$ descends to a well defined
function on $M_\epsilon$. Every point of $M$ has an evenly covered
neighbourhood in $M_\epsilon$, so by shrinking $M_\epsilon$ if
necessary to be the union of these neighbourhoods, we may assume
that any point $x$ of $M_\epsilon$ has an evenly covered
simply connected neighbourhood
$U$ which intersects $M$. Let
$\pi^{-1} U = \cup_k \, \widetilde{U}_k$,
a union of disjoint open sets with biholomorphisms 
$\pi_k: \widetilde{U}_k \to U$.  Now $\hat{F}_+$ restricted to
$\pi^{-1} M$ does descend
to a well defined function on $M$, namely the function $F_+$ which
we started with, and its holomorphic extension to any of the 
biholomorphic neighbourhoods $U_k$ is unique. Hence
 $\hat{F}_+ \circ \pi_k^{-1}$ must be the same function for any $k$,
in other words $\hat{F}_+ \circ \pi^{-1}$ is a well defined function
on $M_\epsilon$.
\end{proof}

%
%
%****************************************************************
%#################################################################
%
\section{Isometric immersions between space forms as pluriharmonic maps} \label{sforms}
To simplify the statement of results, we consider only the three cases
when the target space is a sphere, but note that these contain two of 
the hyperbolic cases, and the other, Case 4 above, can be handled in 
the same way as Case 2.

\begin{theorem} \label{maintheorem}
Let $M$ be a simply connected, paracompact real analytic manifold of 
dimension $m$, and with fixed base point $p$. We can assume
that $M$ is embedded as a totally real submanifold of some
$m$-dimensional complex manifold $M_\cc$.
Let $f: M \to S^{m+k}$ be an immersion with flat normal bundle 
and induced constant
sectional curvature $c$ in one of the following unions of intervals:
\beqas
I_1 =  (-\infty,0),\\
 I_2 =  (0,1), \\
I_3 = (1,\infty),
\eeqas
and $f^\lambda$ the corresponding extended family obtained
by the scheme in \cite{feruspedit1996}.
Then there exists an open submanifold $M_\epsilon$ of $M_\cc$,
containing $M$, and a unique extended family of pluriharmonic maps
into a symmetric space,
$\hat{f}^\lambda: M_\epsilon \to U_i/K_i$, such that the restriction
$\hat{f}^\lambda|_M$ is $f^\lambda$.
The symmetric spaces $U_i/K_i$ corresponding to
$c \in I_i$ are as follows:
\beqas
U_1/K_1 = \frac{SO(m+k,1)}{SO(m) \times SO(k,1)},\\
U_2/K_2 = \frac{SO(m+1,k)}{SO(m) \times SO(k,1)},\\
U_3/K_3 =  \frac{SO(m+k+1)}{SO(m) \times SO(k+1)}.
\eeqas
\end{theorem}

 Theorem \ref{maintheorem} follows from Proposition \ref{theorem2}
below and its analogues for Cases 1 and 3.

Let us consider Case 2 first, a family of isometric immersions 
$M \to S^{m+k+1}$ given by an element of 
\beq \label{spaceform2}
\frac{\hh_{\mu \rho_2}}{\hh_{\mu \rho_2}^0}(M)_{-1}^1,
\eeq
as defined in Section \ref{sketch}.
Observe that (\ref{spaceform2}) looks
rather like 
a special case of the set on the  left hand side of
(\ref{realpluri}), except for the fact that the involution
$\mu$ is not a reality condition (it is not conjugate linear).
However, one can describe this loop group in another way:
replace $\mu$ with the involution ${\tau_2}$ defined by
\bdm
{\tau_2} X(\lambda) := \textup{Ad}_Q (\overline{X(1/\bar{\lambda})}).
\edm
Then ${\tau_2}$ is of the form (\ref{tau}) for the reality condition
$\tau_2^0 (X) := \textup{Ad}_Q(\overline{X})$, and due to the reality condition
${\rho_2}$, which says that $\overline{X(1/\bar{\lambda})} = X(1/\lambda)$ for
$X \in \hh_{{\rho_2}}$, we see that
\bdm
\frac{\hh_{{\rho_2} {\tau_2}}}{\hh_{{\rho_2} {\tau_2}}^0}(M)_{-1}^1 = 
\frac{\hh_{{\rho_2} \mu}}{\hh_{ {\rho_2} \mu}^0} (M)_{-1}^1.
\edm
%**************
\begin{definition} \label{extendim}
An \emph{extended isometric immersion $M \to S^{m+k}$ with induced
constant sectional curvature $c \in (0,1]$} is an element of
$\frac{\hh_{{\rho_2}  {\tau_2}}}{\hh_{{\rho_2}  {\tau_2}}^0}(M)_{-1}^1$, 
where $G = SO(m+k+1, {\bf C})$ and ${\rho_2}$, $\sigma$ and ${\tau_2}$ 
are as defined in this section.
\end{definition}
%*************

It is now easy to prove the following result:
%****************************
\begin{proposition} \label{theorem2}
Let $F$ be a real analytic element of 
$\frac{\hh_{{\rho_2}  {\tau_2}}}{\hh_{{\rho_2}  {\tau_2}}^0}(M)_{-1}^1$, as
given in Definition  \ref{extendim}, where $M$ is a real analytic 
manifold.
Then  there exists a complex manifold $M_\epsilon$ of 
dimension $m$, containing $M$ as a totally real submanifold, 
and a unique extended pluriharmonic map $\hat{F} \in 
\mathcal{P}\frac{\hh_{{\tau_2}}}{\hh_{{\tau_2}}^0}(M_\epsilon)_{-1}^1$, such that
 $\hat{F}|_{M} = F$.
\end{proposition}
\begin{proof}
The big cell  $B\hh$ is a neighbourhood of the identity, so if $p$ is the point in
$M$ at which $F$ is normalised, then there is a neighbourhood $U_p$ of $p$ such that
$F(x)$ takes values in $B\hh$ for all $x$ in $U_p$.  Since the normalisation point is
relevant here, we will use the notation $\hh (M,x)$ to denote maps into $\hh$ which are
normalised at $x$.
Now on $U_p$ we can use the DPW correspondence
(\ref{taucor}) to associate to $F$ a unique extended curved flat 
$F_+ \in \hh_{\rho_2} (U_p, p)_1^1$.  The map $F_+$ is analytic on $U_p$, because
the Birkhoff splitting is analytic.  Theorem \ref{theorem1}
then gives the extension to $\hat{F} \in \mathcal{P} 
\frac{\hh_{{\tau_2}}}{\hh_{{\tau_2}}^0}(U_p ^\epsilon, p)_{-1}^1$, where $U_p \subset U_p ^\epsilon$
and $U_p ^\epsilon$ is open in $M_\cc$. 

For a point $q \in M \setminus U_p$, we consider instead the map
\bdm
R(x,\lambda) := F(q,\lambda)^{-1} F(x,\lambda).
\edm
Now $R^{-1} \dd R = F^{-1}\dd F$,  so $R$ is also an element of 
$\frac{\hh_{{\rho_2}  {\tau_2}}}{\hh_{{\rho_2}  {\tau_2}}^0}(M, q)_{-1}^1$, where here
the normalisation is at $q$ rather than $p$.  We can therefore apply the same argument
to extend $R$ to $\hat{R} \in \mathcal{P}\frac{\hh_{{\tau_2}}}{\hh_{{\tau_2}}^0}(U_q ^\epsilon ,q)_{-1}^1$,
for some neighbourhood $U_q ^\epsilon$ of $q$ in $M_\cc$.  On $U_q ^\epsilon$ we then
define
\bdm
\hat{F}_q(z,\lambda) := F(q,\lambda) \hat{R}(z,\lambda).
\edm
Now $\hat{F}_q$ and $\hat{R}$ have the same Maurer-Cartan form, and therefore, since a map being pluriharmonic is characterised by
 its Maurer-Cartan form, $\hat{F}_q$ is pluriharmonic. In fact $\hat{F}_q$ is an element of 
$\mathcal{P}\frac{\hh_{  {\tau_2}}}{\hh_{{\tau_2}}^0}(U_q ^\epsilon ,p )_{-1}^1$,
where here the ''normalisation'' at $p$ still makes sense for $\hat{F}_q$,
 even if $p$ is not in  $U_q ^\epsilon$,
because $\hat{F}_q$ 
clearly  agrees with $F$ on $M \cap U_q$, and therefore can be extended along $M$ to $p$.

We now want  to check  that 
$\hat{F}_q$ and $\hat{F}$ agree at any point $w \in U_p ^\epsilon \cap U_q ^\epsilon$.
Let $r$ be a point in $U_p \cap U_q$.  It is enough to show that 
\bdm
\hat{F}^{-1}(r,\lambda) \hat{F}_q(w,\lambda) = \hat{F}^{-1}(r,\lambda) \hat{F}(w,\lambda).
\edm
Now both the left and right hand side of this equation are normalised at $r$, because
$r \in M$ and therefore $\hat{F}^{-1}(r,\lambda) = \hat{F}_q^{-1}(r,\lambda)$. 
In fact they are both elements of 
$\mathcal{P} \frac{\hh_{{\tau_2}}}{\hh_{{\tau_2}}^0}(U_q ^\epsilon \cap  U_p ^\epsilon  ,r )_{-1}^1$, which agree with $F^{-1} (r,\lambda) F(x,\lambda)$ for
 $x \in M \cap U_q ^\epsilon \cap  U_p ^\epsilon$.  
 Applying the DPW correspondence (\ref{taucor}) to either one of them
 we get on the right hand side of (\ref{taucor}) the holomorphic curved flat
 $\hat{H}_+ \in \hh(U_q ^\epsilon \cap  U_p ^\epsilon  ,r)_1^1$,
 which has to be the unique holomorphic extension of the curved flat corresponding
 to $F^{-1} (r,\lambda) F(x,\lambda)$.  Hence, by the uniqueness of the correspondence
 (\ref{taucor}), $\hat{F}^{-1}(r,\lambda) \hat{F}_q(w,\lambda)$ and  $\hat{F}^{-1}(r,\lambda) \hat{F}(w,\lambda)$ must be equal.

Repeating the same procedure for any other point $s \in M \setminus \{U_p \cup U_q \}$, 
an identical argument also shows that $\hat{F}_s$ agrees with both $\hat{F}$ and 
$\hat{F}_q$ on the intersections of their respective domains of definition, and thus,
taking $M_\epsilon$ to be the union over all $x \in M$ of the sets $U_x^\epsilon$
we have the required global extension $\hat{F}$.
\end{proof}
Proposition \ref{theorem2} says that a real analytic element of
$\frac{\hh_{{\rho_2}  {\tau_2}}}{\hh_{{\rho_2}  {\tau_2}}^0}(M)_{-1}^1$,
can be extended to a family of pluriharmonic maps into 
$G_{\tau_2^0}/G_{\tau_2^0 \sigma^0}$. Let us identify the
real form $G_{\tau_2^0}$:

%*************
\begin{lemma}  \label{identification}
The real form  $G_{\tau_2^0}$ is isomorphic to $SO(m+1,k,\real)$.
\end{lemma}
\begin{proof} 
Let $J:= \textup{diag}(I_{m+1}, -I_{k})$, and take $SO(m+1,k,\real)$
to be the set of matrices in $GL(m+k+1,\real)$ which satisfy
\bdm
X^tJX = J.
\edm
Let $T := \textup{diag}(I_{m+1}, iI_{k})$.
We  check that $Ad_T$ takes $G_{\tau_2^0}$ isomorphically to
$SO(m+1,k,\real)$: if $X \in G_{\tau_2^0} = SO(m+k+1,\cc)_{\tau_2^0}$,
then the $SO(m+k+1)$ condition $X^tX = I$ implies
\bdm
(Ad_T X)^t J (Ad_T X) = J 
 = T^{-1} X^t X T^{-1} = T^{-1}T^{-1} = J,
\edm
and the condition 
$\tau_2^0 X = \overline{\textup{Ad}_QX} = X$
 implies 
\bdm
\overline{\Ad_R X} = \overline{\Ad_Q \Ad_T \Ad_Q X} = 
\Ad_Q \Ad_{\bar{T}} X = \Ad_T X.
\edm
Thus $\Ad_T X$ preserves $J$ and is real. The converse is similar.
\end{proof}

Now  $G_{\sigma^0}$ is the subgroup of matrices made up of diagonal
blocks which are $m \times m$ and $(k+1)\times (k+1)$ respectively,
and one sees that
\bdm
G_{\tau_2^0 \sigma^0} = SO(m,\real) \times SO(k,1,\real),
\edm
which completes the proof of the first part of Theorem \ref{maintheorem}.

The other two cases can be treated similarly:\\
\noindent \textbf{Case 1}\\
\noindent These are given by elements of 
\bdm
\frac{\hh_{\rho_1 \mu}}{\hh_{\rho_1 \mu}^0}  (M)_{-1}^1.
\edm
 Since $\rho_1$, given by 
\bdm
\rho_1 X(\lambda) := \overline{X(-\bar{\lambda})},
\edm
 is a reality condition
along $i \real$ rather than $\real$ we cannot define ${\tau}$ as we did
in Case 2. 
Instead we define 
\bdm
{\tau}_1 X(\lambda) := \textup{Ad}_Q \textup{Ad}_P \overline{X(1/\bar{\lambda})},
\edm
where we recall that $(\sigma X)(\lambda) := \textup{Ad}_PX(-\lambda)$,
and $(\mu X)(\lambda) := \textup{Ad}_Q(X(1/\lambda))$.

Again, $\tau_1$ is an involution of the type (\ref{tau}) for the
reality condition 
\bdm
\tau^0_1 X = \textup{Ad}_Q \textup{Ad}_P \overline{X},
\edm
and, for $X \in \hh_{\rho_1}$ we have
\beqas
\tau_1(X(\lambda)) &=& \textup{Ad}_Q \sigma \rho_1 X(1/\lambda)\\
&=& \textup{Ad}_Q X(1/\lambda)\\
&=& \mu X(\lambda).
\eeqas
From this it follows that
\bdm
\hh_{\rho_1 \tau_1}= \hh_{\rho_1 \mu},
\edm
and we can proceed as in Case 2. Note that, although $\rho_1$
is not a reality condition along $\real$, as $\rho$ was, the only thing
that mattered in our construction of pluriharmonic maps  from 
curved flats was that $\rho$ should be a positive  reality condition,
and this is the case for $\rho_1$.

Evidently, elements of 
$\frac{\hh_{\rho_1 \mu}}{\hh_{\rho_1 \mu}^0}  (M)_{-1}^1$
extend to pluriharmonic maps into $G_{\tau_1^0}/G_{\tau_1^0 \sigma^0}$.
Now $\textup{Ad}_Q \textup{Ad}_P = \textup{Ad}_{QP}$, where
$QP = \textup{diag}(I_{m+k,-1})$, and so the same argument given in
Lemma \ref{identification} shows that the real form $G_{\tau_1^0}$
is isomorphic to $SO(m+k,1,\real)$, while the subgroup 
$G_{\tau_1^0 \sigma^0}$ is $SO(m,\real) \times SO(k,1,\real)$.

\noindent \textbf{Case 3}\\
\noindent
Here the relevant loop group is $\hh_{\rho_3 \mu}$,
and $\rho_3$ is the negative  reality condition:
\bdm
(\rho_3 X)(\lambda) := \overline{X(1/\bar{\lambda})}.
\edm
  In this case
we define a positive involution $\hat{\rho}_3$ by
\beqas
\hat{\rho}_3X(\lambda) &:=& \rho_3 \mu X(\lambda)\\
&=& \textup{Ad}_Q \overline{X(\bar{\lambda})},\\
\eeqas
and one easily verifies that
\bdm
\hh_{\hat{\rho}_3 \mu}= \hh_{\rho_3 \mu}.
\edm
This is then the same as the first case, as $\hat{\rho}_3$ is an
involution of the form (\ref{rho}) for the reality condition
\bdm
\hat{\rho}^0_3 X = \textup{Ad}_Q \overline{X}.
\edm
 Proceeding as in  Case 2, we replace $\mu$ with $\tau_3$ defined by:
\beqas
({\tau_3}X)(\lambda) &:=& (\mu \hat{\rho}_3X)(\lambda) \\
&=& \overline{X(1/\bar{\lambda})}.
\eeqas
Thus an element of $\hh_{{\tau_3}}(M_\epsilon)_{-1}^1$ is 
a pluriharmonic map into 
\bdm
U/K = SO(m+k+1,\real)/(SO(m,\real) \times
SO(k+1,\real),
\edm
 and this completes the proof of Theorem \ref{maintheorem}.

\subsection{Some remarks on Theorem \ref{maintheorem}}
\begin{enumerate}
\item
A similar result  
holds even if $M$ is not simply connected, the only
difference being that the associated family $f^\lambda$ is in general
only defined on the universal cover.
\item
The original map $f$ is actually a map into the
last of these symmetric spaces, $U_3/K_3$, which is not
isomorphic to the first two, and so the result might seem
odd, but the explanation is that the pluriharmonic maps
are obtained by evaluating the extended map at $\lambda \in S^1$, while
in the first two cases, the isometric immersions are obtained 
from values  of $\lambda$ in $i\real^*$ and $\real^*$ respectively.
On the intersections, $\lambda = \pm 1$ and $\lambda = \pm i$, 
the maps take their values in the intersections of the relevant 
symmetric spaces.
\end{enumerate}

\subsection{Totally geodesic immersions}  \label{application}
Recall from Section \ref{part2} that totally geodesic immersions into
the sphere and hyperbolic space were obtained from the extended
families of immersions evaluated at $\lambda = \pm 1$,
but that given a totally geodesic immersion
$f$ into the sphere, 
 there is no canonical way to insert the
parameter $\lambda$  to obtain the extended family.
 Nor is there an obvious way to tell
whether such a family exists for $f$ (in contrast to other values
of the induced curvature $c$).
To see that the question is meaningful,
consider a totally geodesic embedding of $S^n$ into $S^{2n-1}$.  This cannot belong
to one of the families of Case 2 above, that is with constant curvature
varying in the interval $(0,1]$, because the immersion condition is preserved within
these families, and, as noted before,
 there is no global isometric
immersion of a sphere  $S_c^n$, of radius $1/\sqrt{c}$ for
$0<c<1$, into $S^{2n-1}$ \cite{moore1}.

As an application of Theorem \ref{maintheorem}, we characterise
in terms of pluriharmonic maps, which real analytic totally
geodesic submanifolds of $S^{m+k}$ can be extended to one of
the families discussed here. The analogue holds for totally
geodesic submanifolds of $H^{m+k}$.
The point here is that, while there is no canonical way to 
insert the parameter $\lambda$ into a totally geodesic
immersion, there \emph{is} for a pluriharmonic map.

Let $M$ be as in Theorem \ref{maintheorem}, suppose 
$f: M \to S^{m+k}$ is a totally geodesic immersion with flat normal bundle,
and regard $f$ as a map into
 $SO(m+k+1)/SO(m) \times SO(k)$, given by its equivalence class of
adapted frames.
  Since the second fundamental form of 
$f$ is zero, one may check that $f$ satisfies both 
the reality conditions $\tau_2^0$ and $\tau_3^0$, which,
combined with $\sigma^0$ define the symmetric spaces
 $SO(m+1,k)/(SO(m) \times SO(k,1))$ and
$SO(m+k+1)/(SO(m) \times SO(k+1))$ respectively.  
Around a point
$p \in M \subset M_\cc$, the real analytic function $f$ can be extended
locally to a 
 pluriharmonic map from $M_\cc$ into either of these symmetric
spaces, since  such maps  are locally just
the real parts of holomorphic functions.  These extensions are not unique,
however, since the notion of real part depends on your choice of coordinates
for the target space. 

The following is an immediate corollary of Theorem \ref{maintheorem}.
%***************************************
\begin{proposition} \label{corollary1}
If $f: M \to S^{m+k}$ is a real analytic totally geodesic immersion with flat normal bundle, then:
\begin{enumerate}
\item
$f = f^1$ for some  extended isometric 
immersion $f^\lambda$ with 
$c_\lambda \in (0,1]$ if and only if it has an extension to a pluriharmonic map 
$\hat{f} : M_\cc \to SO(m+1,k)/(SO(m) \times SO(k,1))$ such that, if
$\hat{f}^\lambda$ is the corresponding extended pluriharmonic map,
 then $\hat{f}^\lambda|_M$ is 
fixed by the involution $\rho_2$.
\item
$f = f^1$ for some extended isometric immersion
 $f^\lambda$ with 
$c_\lambda \in [1,\infty)$ if and only if 
it has an extension to a pluriharmonic map 
$\hat{f} : M_\cc \to SO(m+k+1)/(SO(m) \times SO(k+1))$
such that, if
$\hat{f}^\lambda$ is the corresponding extended pluriharmonic map,
 then $\hat{f}^\lambda|_M$ is 
fixed by the involution $\hat{\rho}_3$.
\end{enumerate}
\end{proposition}
\vspace{.5cm}
\section*{Acknowledgements}
This research was partially 
supported by DFG grant DO 776. The author would like to 
thank Josef Dorfmeister for encouraging him
to work on this project and for many helpful discussions.
Thanks are also due  to  Jost-Hinrich Eschenburg and Wayne
Rossman for helpful discussions and feedback.

\bibliographystyle{abbrv} 

\bibliography{mybib}

\end{document}